\documentclass[a4paper]{amsart}
\usepackage{amsmath, amsfonts,amsthm,amssymb,amscd, verbatim,graphicx,color,multirow,booktabs, caption, bm,chngcntr,adjustbox,longtable,mathtools,microtype, xfrac}
\usepackage[top = 78pt, bottom = 72pt, inner=108pt, outer=108pt]{geometry}
\usepackage{hyperref}
\usepackage[dvipsnames, table]{xcolor}
\hypersetup{colorlinks=true,linkcolor=Maroon, citecolor=OliveGreen}

\usepackage[capitalise,noabbrev]{cleveref}
\usepackage[shortlabels]{enumitem}

\newcommand{\Z}{\mathbb{Z}}

\newtheorem{theorem}{Theorem}[section]
\newtheorem{lemma}[theorem]{Lemma}

\newtheorem{proposition}[theorem]{Proposition}

\title[On generalisations of conciseness]
 {On generalisations of conciseness}

\author[A.~Zozaya]{Andoni Zozaya} 
\address{Department of Statistics, Computer Science and Mathematics, Public University of Navarra (UPNA) \& Institute for Advanced Materials and Mathematics (INAMAT$^2$), \linebreak Arrosadia Campus, 31006 Pamplona, Spain}
\email{andoni.zozaya@unavarra.es}

\begin{document}

\begin{abstract}
Based on the notions of conciseness and semiconciseness, we show that these properties are not equivalent by proving that a word originally presented by Ol'shanskii is semiconcise but not concise. We further establish that every \(1/m\)-concise word is semiconcise by proving that when the group-word $w$ takes finitely many values in $G$, the iterated commutator subgroup \([w(G), G, \stackrel{(m)}{\dots}, G]\) is finite for some \(m \in \mathbb{N}\) if and only if \([w(G), G]\) is finite.

\medskip
\noindent{\scshape MSC 2020.} 20F10

\noindent{\scshape Keywords.} Group word -- Verbal subgroup -- Conciseness -- Semiconciseness
\end{abstract}


\thanks{The author is supported by the project PID2020-117281GB-I00 (Spanish Government, partially funded with ERDF) and by the research group 244 \textit{\'Algebra. Aplicaciones} (Public University of Navarra).}

\maketitle 

\section{Introduction}

A \emph{word} in $k$ variables is an element $w = w(x_1, \dots, x_k)$ of the free group on $k$ generators $F(x_1, \dots, x_k)$. Given a group $G$, such a word induces a natural \emph{word map} $w\colon G^k \rightarrow G$ through substitution. As defined by P.~Hall, the word $w$ is said to be \emph{concise} in a class of groups $\mathcal{C}$ if, for every $G \in \mathcal{C}$, the finiteness of the set of \emph{$w$-values}
\[
w\{G\}  = \{ w(g_1, \dots, g_k) \mid g_i \in G \}
\]
implies that the verbal subgroup $w(G) = \langle w\{G\} \rangle$ is also finite.

Although this property was originally conjectured to hold for all words in the class consisting of all groups, counterexamples were constructed first by Ivanov~\cite{Ivanov}, and later by Ol’shanskii~\cite{Old}. Nonetheless, conciseness, and especially the question of whether every word is concise in the class of residually finite groups, remains actively studied.

A weaker variant, namely \emph{semiconciseness}, was introduced in~\cite{DST}. A word $w$ is semiconcise in a class of groups $\mathcal{C}$ if, for every $G \in \mathcal{C}$, the finiteness of $w\{G\}$ implies that the commutator subgroup $[w(G), G]$ is finite.  For instance, \emph{Engel words}, the two-variant words defined recursively by $[x{,}_{1} y ] =[x,y] $ and
\[[x,_{\, n} y] = \left[[x,_{\, n-1} y], \, y \right], \quad \, n \geq 3,\]
are semiconcise in the class of all groups; see~\cite{FMT}. Yet for $n \geq 5$ it is still unknown whether they are concise in the class of all groups (in the class of residually finite groups, all Engel words have been proved concise; see~\cite{DMS}).

In Section~\ref{section semiconcise}, we show that conciseness and semiconciseness are, indeed, distinct properties. Specifically, we prove that the word that appears in Ol’shanskii’s aforementioned counterexample, referred to here as \emph{Olshanskii's word}, is semiconcise (in the class of all groups), despite not being concise. Thereby we answer a question posed in~\cite{DST}.

\begin{theorem}
\label{thm 1}
Let $w_o$ be Ol'shanskii's word. Suppose that $w_o\{G\}$ is finite for a group $G$. Then $[w_o(G), G]$ is finite.
\end{theorem}

A further generalisation was proposed in~\cite{DGM}. For each integer $m \geq 1$, a word $w$ is called \emph{$1/m$-concise} in $\mathcal{C}$ if, for every $G \in \mathcal{C}$, the finiteness of $w\{G\}$ implies that the iterated commutator subgroup $[w(G), G, \stackrel{(m)}{\dots}, G]$ is finite. This gives rise to a purported hierarchy of properties, where $1/m$-conciseness implies $1/n$-conciseness for all $n \geq m$. In Section~\ref{section m concise}, we prove that semiconciseness is equivalent to $1/m$-conciseness for all $m$.

\begin{theorem}
\label{thm 2}
A word $w$ is semiconcise in a group $G$ if and only if it is $1/m$-concise in $G$ for every $m \geq 1$.
\end{theorem}


\section{Proof of Theorem~\ref{thm 1}}
\label{section semiconcise}


In this section, we consider the two-variant group-word $w_o$ originally constructed by Ol'shanskii as an example of a non-concise word; see~\cite[Theorem 39.7]{Old}.  This specific word has since been further studied, and it has been shown that it is concise when restricted to the class of residually finite groups; see~\cite[Theorem 3.1]{PS}.

The word \( w_o \) is defined using an auxiliary word: let \( d \geq 1 \) be an integer and define
\[
v(x, y) = \left[ \left[x^d, y^d \right]^d, \left[y^d, x^{-d} \right]^d \right].
\]
The word \( w_o(x, y) \) is constructed as follows:
\[
w_o(x, y) = [x, y] \cdot v(x, y)^n \cdot [x, y]^{\varepsilon_1} \cdot v(x, y)^{n+1} \cdots [x, y]^{\varepsilon_{h-1}} \cdot v(x, y)^{n+h-1},
\]
where the exponents \( \varepsilon_i \) follow the periodic pattern:
\[
\varepsilon_{10k+1} = \varepsilon_{10k+2} = \varepsilon_{10k+3} = \varepsilon_{10k+5} = \varepsilon_{10k+6} = 1,
\]
\[
\varepsilon_{10k+4} = \varepsilon_{10k+7} = \varepsilon_{10k+8} = \varepsilon_{10k+9} = \varepsilon_{10k+10} = -1,
\]
for \( k = 0, 1, \dots, \sfrac{h - 1}{10} \), assuming \( h \equiv 1 \pmod{10} \) and \( h > 50000 \). \\

We will use the following well-known results:

\begin{lemma}
\label{lemma: fin conjugacy class}
Let $G$ be a group and let $w$ be a word such that $w\{G\}$ is finite. 
\begin{itemize}
\item[\textup{(i)}] For every $h$ in $w(G)$ the conjugacy class $h^G$ is finite.
\item[\textup{(ii)}] The derived subgroup $w(G)'$ is finite.
\end{itemize}
\end{lemma}
\begin{proof}
(i) Denote by $w\{G\}^{*\ell}$ the elements that can be written as the product of $\ell$ $w$-values or inverses of $w$-values. Then, $h \in w\{G\}^{*\ell}$ for some $\ell$. Note that if $ w(g_1, \dots, g_k)$ is a $w$-value, then 
\[ w(g_1, \dots, g_k)^t  = w(g_1^{t}, \dots, g_k^{t}) \]
is again a $w$-value for any $t \in G$. Thus, $ h^G \subseteq w\{G\}^{*\ell}$, and $|h^G| \leq |w\{G\}|^\ell$ is finite. \\ 
\linebreak
(ii) Since $|G : C_G(h)| = |h^G|$ is finite for every $w$-value, then
\[C_G(w(G)) = \bigcap_{h \in w\{G\}} C_G(h)\]
has also finite index in $G$. Thus, $|w(G) : Z(w(G))|$ is finite, and by Schur's lemma,  $w(G)'$ is also finite; see~\cite[Theorem 10.1.4]{Robinson}.
\end{proof}

From the construction of $w_o$, it follows that when the auxiliary word $v$ is a law in $G$, \textit{e.g.} when $G$ is metabelian, then $w_o$ coincides with the commutator word $\gamma_2$.

\begin{lemma}
\label{lemma: metabelian}
Let $G$ be a metabelian group. Then $w_o\{G\} = \gamma_2\{G\}$.
\end{lemma}

\begin{proof}[proof of Theorem \ref{thm 1}]
Let $G$ be a group such that $w_o\{G\}$ is finite. In view of Lemma~\ref{lemma: fin conjugacy class}, we may assume that $w_o(G)$
torsion-free and abelian, after quotienting by suitable finite subgroups.

We shall prove that $[w, g] = 1$ for every $w \in w_o\{G\}$ and $g \in G$. Fix $g \in G$, and consider the metabelian group $M = w_o(G) \rtimes \langle g \rangle \leq G$. Since $w_o\{M\} \subseteq w_o\{G\}$ is finite, and every word is concise in the class of virtually metabelian groups (\textit{cf.}~\cite[Corollary 2]{TS}), it follows that $w_o(M)$ is finite. But $w_o(G)$ is torsion-free, and so $w_o(M) = \{1\}$.

Finally, by Lemma~\ref{lemma: metabelian}, we have that
\[
[w, g] \in \gamma_2(M) = w_o(M) = \{1\},
\]
for every $w \in w\{ G \}$, as claimed.
\end{proof}


\section{Proof of Theorem~\ref{thm 2}}
\label{section m concise}


Theorem~\ref{thm 2} follows immediately from this result:

\begin{proposition}
Let $G$ be a group, and let $w$ be a word. Suppose that $w\{G\}$ is finite and that $[w(G), G, \stackrel{(m)}{\dots}, G]$ is finite for some integer $m\geq 1$. Then, $[w(G), G]$ is finite. 
\end{proposition}
\begin{proof}
Denote $W_i = [w(G), G, \stackrel{(i)}{\dots}, G]$. We shall prove for any $n \geq 2$ that if $W_n$ is finite, then $W_{n-1}$ is also finite.

\medskip

Note that for every normal subgroup $N \unlhd G$, 
\[ [w(G/N), G/N, \stackrel{(n)}{\dots}, G/N]]  = \frac{[w(G), G, \dots, G] N}{N}\cong \frac{[w(G), G, \dots, G]}{[w(G), G, \dots, G] \cap N}.\]

Therefore, since $W_{n}$ is finite, we may suppose, without loss of generality, that $[w(G), G, \stackrel{(n)}{\dots}, G] = 1$, \textit{i.e.}, $W_{n-1} = [w(G), G, \stackrel{(n-1)}{\dots}, G]$ is central in $G$. Moreover,  by Lemma~\ref{lemma: fin conjugacy class}~(ii), since $w(G)'$ is finite, we may further assume that $w(G)$ is abelian.

\medskip

Let $W_{n-2}$ be generated by say $\{ w_1, \dots, w_d\} \subseteq w(G)$. Any element $g \in W_{n-1}$ can be written as
\begin{equation*}
\begin{split}
g & = [\mu_1(w_1, \dots, w_d), g_1] \cdot \dots \cdot [\mu_\ell(w_1, \dots, w_d), g_\ell] \\
 & = \prod_{j=1}^{\ell} \left[ \prod_{i=1}^d w_i^{k_{ij}}, g_j \right] = \prod_{i=1}^d \prod_{j=1}^{\ell} \left[  w_i, g_j^{k_{ij}} \right] \\ & = \prod_{i=1}^d \left[ w_i, \prod_{j=1}^\ell g_j^{k_{ij}} \right],
\end{split}
\end{equation*}
where $\mu_j$ is a word in $d$-variables, $g_j \in G$ and $k_{ij} \in \Z$. The above equalities follow because $w(G)$ is abelian and the commutators are bilinear, as $W_{n-1}$ is central in $G$.

\medskip 

In particular, there exist elements \( \tilde{g}_i \in G \), for \( i = 1, \dots, d \), such that
\[
g = [w_1, \tilde{g}_1] \cdots [w_d, \tilde{g}_d].
\]

Hence,
\[
W_{n-1} \subseteq \mathcal{C}_1 \cdots \mathcal{C}_d,
\]
where
\[
\mathcal{C}_i := \{ [w_i, g] = w_i^{-1} w_i^g \mid g \in G \} = w_i^{-1} w_i^G.
\]
Since each element of \( w(G) \) has a finite conjugacy class by Lemma~\ref{lemma: fin conjugacy class}~(i), each set \( \mathcal{C}_i \) is finite. Thus, \( W_{n-1} \) is finite as required.
\end{proof}


\bibliographystyle{plain}
\bibliography{biblio}

\begin{thebibliography}{1}

\bibitem{DGM}
C.~Delizia, M.~Gaeta, and C.~Monetta.
\newblock On generalized concise words.
\newblock {\em J. Group Theory}, 28:475--487, 2025.

\bibitem{DST}
C.~Delizia, P.~Shumyatsky, and A.~Tortora.
\newblock On semiconcise words.
\newblock {\em J. Group Theory}, 23:629--639, 2020.

\bibitem{DMS}
E.~Detomi, M.~Morigi, and P.~Shumyatsky.
\newblock Words of {Engel} type are concise in residually finite groups.
\newblock {\em Bull. Math. Sci.}, 9: art. 1950012, 2018.

\bibitem{FMT}
G.~A. Fern\'andez-Alcober, M.~Moriagi, and G.~Traustason.
\newblock A note on conciseness of {Engel} words.
\newblock {\em Comm. Algebra}, 40:2570--2576, 2012.

\bibitem{Ivanov}
S.~V. Ivanov.
\newblock P. {H}all’s conjecture on the finiteness of verbal subgroups.
\newblock {\em Soviet Math. (Iz. VUZ)}, 33:59--70, 1989.

\bibitem{Old}
A.~Yu. Ol'shanskii.
\newblock {\em Geometry of Defining Relations in Groups}, volume~70 of {\em Mathematics and its Applications (Soviet Series)}.
\newblock Kluwer Academic Publishers, Dordrecht, 1991.

\bibitem{PS}
M.~Pintonello and P.~Shumyatsky.
\newblock On conciseness of the word in {O}lshanskii’s example.
\newblock {\em Arch. Math.}, 122:241--247, 2024.

\bibitem{Robinson}
D.~J.~S. Robinson.
\newblock {\em A course in the theory of groups}.
\newblock Springer-Verlag, New York, 1982.

\bibitem{TS}
R.~F. Turner-Smith.
\newblock Finiteness conditions for verbal subgroups.
\newblock {\em J. Lond. Math. Soc.}, 41:166--176, 1966.

\end{thebibliography}

\end{document}